\documentclass[12pt]{amsart}

\usepackage{amssymb}

\usepackage{enumerate}

\usepackage{graphicx}

\makeatletter
\@namedef{subjclassname@2010}{%
  \textup{2010} Mathematics Subject Classification}
\makeatother

\newtheorem{thm}{Theorem}
\newtheorem{cor}[thm]{Corollary}

\theoremstyle{definition}
\newtheorem{defin}{Definition}

\begin{document}

\baselineskip=17pt



\title[Counting rational points]{Counting rational points near planar curves}

\author[A. Gafni]{Ayla Gafni}
\address{109 McAllister Bldg\\ Penn State University\\
University Park, PA 16802}
\email{gafni@math.psu.edu}

\date{}

\begin{abstract}
We find an asymptotic formula for the number of rational points near planar curves.  More precisely, if $f:\mathbb{R}\rightarrow\mathbb{R}$ is a sufficiently smooth function defined on the interval $[\eta,\xi]$, then the number of rational points with denominator no larger than $Q$ that lie within a $\delta$-neighborhood of the graph of $f$ is shown to be asymptotically equivalent to $(\xi-\eta)\delta Q^2$.\end{abstract}



\maketitle

\section{Introduction}

In this paper, we give an explicit asymptotic formula for the number of rational points with bounded denominator near a sufficiently smooth planar curve.  This result expands on Theorem 3 of \cite{VaVe06}, and it may be able to provide quantitative information about Khinchin-type manifolds.  

The results in this paper are motivated by the convergence side of Khinchin theory, and so we will begin with an overview of the relevant points therein.    We say that $\psi : \mathbb{R}^+ \rightarrow \mathbb{R}^+$ is an \emph{approximating function} if it is decreasing and satisfies $\psi(x)\rightarrow 0$ as $x\rightarrow\infty$.  Given an approximating function $\psi$, we say that a point $(y_1, \ldots, y_n) \in \mathbb{R}^n$ is \emph{simultaneously $\psi$-approximable} if there exist infinitely many $q\in\mathbb{N}$ such that 
\begin{equation}\label{approximable if}\max_{1\le i \le n} || qy_i || \le \psi(q).\end{equation}
Here $||x|| = \min_{m\in\mathbb{Z}} |x-m|$.  We denote by $\mathcal{S}(\psi)$ the set of all simultaneously $\psi$-approximable points in $\mathbb{R}^n$.  Khinchin's theorem gives a criterion for the $n$-dimensional Lebesgue measure $|\cdot |_{\mathbb{R}^n}$ of $\mathcal{S}(\psi)$, namely
$$\left|\mathcal{S}(\psi) \right|_{\mathbb{R}^n} = \left\{
\begin{array}{l l} 0 & \text{ if }\sum_{q \ge 1} \psi(q)^n <\infty \\ 
\text{\sc Full} & \text{ if }\sum_{q \ge 1} \psi(q)^n = \infty
\end{array}\right.,$$
where ``{\sc Full}'' means that the complement of the set has measure 0.

Current research in metric Diophantine approximation focuses on expanding this theorem to $m$-dimensional manifolds in $\mathbb{R}^n$.  Let $\mathcal{M}\subset\mathbb{R}^n$ be a manifold and denote the induced Lebesgue measure on $\mathcal{M}$ by $|\cdot |_{\mathcal{M}}$.  We say that $\mathcal{M}$ is of \emph{Khinchin type for convergence} if $|\mathcal{M} \cap \mathcal{S}(\psi)|_{\mathcal{M}} =  0$ for any approximating function $\psi$ with $\sum_{q \ge 1} \psi(q)^n <\infty$.   Similarly, we say that  $\mathcal{M}$ is of \emph{Khinchin type for divergence} if $|\mathcal{M} \cap \mathcal{S}(\psi)|_{\mathcal{M}} =  $ {\sc Full} for any approximating function $\psi$ with $\sum_{q \ge 1} \psi(q)^n = \infty$.

In this paper we are specifically concerned with curves in $\mathbb{R}^2$.  It is established by Beresnevich et al.  in \cite{BeDiVeVa07} that any $C^{(3)}$ non-degenerate planar curve is of Khinchin-type for divergence.  Vaughan and Velani establish in [5] that such curves are also of Khinchin-type for convergence.  The proof of the convergence case relies on an upper bound on the number of rational points near the curve.  The intuition is that if there are not many rational points near the curve, then we cannot have many approximable points.  This paper provides an asymptotic formula for the number of rational points near a curve.  These results may lead to information about the growth of the number of solutions to (\ref{approximable if}) with $q\le Q$, as $Q\rightarrow\infty$.

\section{Statement of results}

\begin{defin} Let $\eta, \xi\in\mathbb{R}, \ \eta<\xi,\  I = [\eta,\xi]$  and $f : I \rightarrow \mathbb{R}$ be such that $f''$ is continuous and bounded away from $0$ on $I$.   For $Q \ge 1$ and $0<\delta<1/2$, define 
$$N(Q,\delta) := \mbox{card}\{ (a,q) \in \mathbb{Z} \times \mathbb{N} : 1\le q \le Q, \eta q < a \le \xi q, ||qf(a/q)|| < \delta\}.$$
\end{defin}

When dealing with rational points in $\mathbb{R}^n$, we consider the ``denominator'' of the point to be the least common denominator of the coordinates of the point.   Then $N(Q,\delta)$ counts the number of rational points within a $\delta$-neighborhood of the curve graphing $f$, where we require that the denominator of the points be no more than $Q$.  When we apply our results to Khinchin theory, the parameter $\delta$ will be replaced by a suitable approximating function $\psi(q)$.  It is therefore reasonable, when finding asymptotic formulae, to bound $\delta$ from below in terms of $Q$.  

The computations are easier when all values of $q$ are of the same order of magnitude, so we will in fact be working with a slightly different object, namely
$$\widetilde{N}(Q,\delta) := \mbox{card}\{ (a,q) \in \mathbb{Z} \times \mathbb{N} : Q< q \le 2Q, \eta q < a \le \xi q, ||qf(a/q)|| < \delta\}.$$
Theorem 1 gives an explicit asymptotic formula for $\widetilde{N}(Q,\delta)$.  We translate this back to $N(Q,\delta)$ in Theorem 2.  
\begin{thm} Suppose that  $0 < \theta < 1$ and $f'' \in \mbox{\emph{Lip}}_\theta([\eta, \xi])$.  If $Q^{-\frac{1+\theta}{3-\theta}+\varepsilon}\le\delta < 1/2$, then 
$$\widetilde{N}(Q,\delta) = 3(\xi - \eta) \delta Q^2 + E(Q,\delta),$$
where the error term satisfies\footnote{The first range of $\delta$ will not occur when $\theta\le1/2$.}
\begin{equation}\label{official error}E(Q,\delta) \ll  \left\{ \begin{array}{l c l}
\delta^{2/3}Q^{5/3}(\log Q)^{2/3} & \mbox{ if } & \delta>>Q^{\frac{1-2\theta}{2-\theta}}(\log Q)^{-\frac{5-\theta}{2-\theta}} \\ 
\delta^{\frac{2}{5-\theta}}Q^{\frac{3(3-\theta)}{5-\theta}}& \mbox{ if } & \delta\ll Q^{\frac{1-2\theta}{2-\theta}}(\log Q)^{-\frac{5-\theta}{2-\theta}}.
\end{array}\right.\end{equation}
\end{thm}

\begin{thm}  For $\theta, f,$ and $\delta$ as above, we have 
$$N(Q,\delta) =(\xi - \eta) \delta Q^2 + F(Q,\delta),$$
where $F(Q,\delta)$ satisfies the bound given by (\ref{official error}).
\end{thm}
\begin{cor} For $\theta, f,$ and $\delta$ as above, we have
$$\widetilde{N}(Q,\delta) \sim 3(\xi - \eta) \delta Q^2.$$
\end{cor}
\begin{cor}
For $\theta, f,$ and $\delta$ as above, we have
$$N(Q,\delta) \sim (\xi - \eta) \delta Q^2.$$
\end{cor}

\section{Proof of Theorem 1}  

For convenience we extend the definition of $f$ to $\mathbb{R}$ by defining $f(\beta)$ to be $\frac{1}{2} (\beta - \xi)^2 f''(\xi) + (\beta -\xi) f'(\xi) + f(\xi)$ when $\beta > \xi$ and $\frac{1}{2} (\beta - \eta)^2 f''(\xi) + (\beta - \eta) f'(\xi) + f(\xi)$ when $\beta < \eta$.  Note that then $f''\in \mbox{Lip}_\theta (\mathbb{R})$ and $f''$ is still bounded away from 0 and is bounded.  

We follow the methods of the proof of Theorem 3 in \cite{VaVe06}.   Let $K$ be a sufficiently large integer that will be determined later.  Let $S_K^+(\alpha), S_K^-(\alpha)$ be the Selberg functions for the interval $J = (-\delta,\delta)$.  These functions are trigonometric polynomials of degree at most K with the properties that $S_K^-(\alpha) \le \chi_{J}(\alpha) \le S_K^+(\alpha)$ for all $\alpha$ and $\int_{\mathbb{T}} S_K^{\pm}(\alpha)\, d\alpha = 2\delta \pm \frac{1}{K+1}$.  See Section 7.2 of \cite{Montgomery94} for more details about these functions.

From the definition of $\widetilde N(Q,\delta)$ and the properties of the Selberg functions, we see that 
\begin{align*}\widetilde{N}(Q,\delta) & = \sum_{Q < q \le 2Q} \sum_{ \eta q < a \le \xi q} \chi_{J} (||qf(a/q)||) \\  & 
\le  \sum_{Q< q \le 2Q} \sum_{ \eta q < a \le \xi q} S_{K}^+ (qf(a/q)) \\ &
=    \sum_{Q< q \le 2Q} \sum_{ \eta q < a \le \xi q} \sum_{k=-K}^K \widehat{S}_{K}^+ (k) e(kqf(a/q)) \\ & = N_0^{+} +  N_1^{+},
\end{align*}
where 
$$ N_0^{+} := \sum_{Q< q \le 2Q} \sum_{ \eta q < a \le \xi q}\widehat{S}_{K}^+ (0)$$ and
$$ N_1^{+} := \sum_{0 < |k| \le K} \widehat{S}_{K}^+ (k)\sum_{Q< q \le 2Q} \sum_{ \eta q < a \le \xi q}  e(kqf(a/q)).$$

\noindent We wish to find a suitable upper bound for $N_0^+$.  Recall that $\widehat{S}_{K}^+ (0) = \int_{\mathbb{T}} S_K^{+}(\alpha)\,d\alpha = 2\delta + \frac{1}{K+1} $.  Since there are at most $(\xi - \eta)q+ 1$ integers in the interval $\left(\eta q, \xi q\right]$,  we have 
\begin{align*} 
N_0^{+}  & \le  \left(2\delta + \frac{1}{K+1} \right) \left((\xi-\eta)\frac{Q(3Q+1)}{2} + Q\right) \\ &
= 3(\xi - \eta) \delta Q^2 + (\xi - \eta + 2)\delta Q + \frac{3(\xi - \eta)Q^2}{2(K+1)} + \frac{(\xi - \eta +2)Q}{2(K+1)}  \\ &
= 3(\xi - \eta) \delta Q^2 + O\left(\delta Q + K^{-1}Q^2 \right).
\end{align*}

Using $S_K^-$ in place of $S_K^+$, we similarly find that 
$$\widetilde{N}(Q,\delta) \ge N_0^{-}  + N_1^{-} $$

\noindent where 
$$ N_0^{-} := \sum_{Q< q \le 2Q} \sum_{ \eta q < a \le \xi q}\widehat{S}_{K}^- (0) \ge  3(\xi - \eta) \delta Q^2 + O\left(\delta Q + K^{-1}Q^2 \right)$$
 and
$$ N_1^{-} := \sum_{0 < |k| \le K} \widehat{S}_{K}^- (k)\sum_{Q< q \le 2Q} \sum_{ \eta q < a \le \xi q}  e(kqf(a/q)).$$

It can easily be shown that $|\widehat{S}_{K}^\pm (k)| \le |\widehat{S}_{K}^\pm (0)| \ll  \delta + K^{-1}$.  For convenience we define
$$N_1 := \sum_{0 < |k| \le K}(\delta + K^{-1})\left|\sum_{Q< q \le 2Q} \sum_{ \eta q < a \le \xi q}  e(kqf(a/q))\right|.$$
It then follows that $N_1^+, N_1^- \ll  N_1$.  Thus from the above analysis, we see that
$$\widetilde{N}(Q,\delta) = 3(\xi - \eta) \delta Q^2 +  O\left(N_1 + \delta Q + K^{-1}Q^2 \right).$$

\noindent In other words,
\begin{equation}\label{initial error}E(Q,\delta) := \widetilde{N}(Q,\delta) - 3(\xi - \eta) \delta Q^2 \ll  N_1 + \delta Q + K^{-1} Q^2.\end{equation}
In order to find an upper bound for $E(Q,\delta)$,  we need to compute an upper bound for $N_1$ in terms of $\delta, K,$ and $Q$.  This part of the proof is entirely similar to the proof of Theorem 3 in  \cite{VaVe06}, and many of the details are omitted here.  

Consider the function $F(\alpha) = kqf(a/q)$, which has derivative $kf'(a/q)$.  Given $k$ with $0 < |k| \le K$, we define
\begin{align*}
H_- & = \lfloor \inf kf'(\beta)\rfloor - 1, \ \ \  H_+ = \lceil\sup kf'(\beta)\rceil+ 1, \\
h_- & = \lceil\inf kf'(\beta)\rceil+ 1, \ \ \ \ h_+ = \lfloor \sup kf'(\beta)\rfloor - 1
\end{align*}
where the extrema are taken over the interval $[\eta, \xi]$.  By Lemma 4.2 of \cite{Vaughan97}, we have that 
$$\sum_{ \eta q < a \le \xi q}  e(kqf(a/q)) = \sum_{H_- \le h \le H_+} \int_{\eta q}^{\xi q} e(kqf(\alpha/q) - h\alpha)\, d\alpha + O(\log(2+H)),$$
where $H = \max(|H_-|,  |H_+|)$.  So we have,
\begin{equation}N_1 = N_2 + O\left(\sum_{0 < |k| \le K} (\delta + K^{-1}) \sum_{Q< q \le 2Q} \log(2+H)\right),
\label{N1 to N2}\end{equation}
where $$N_2 =   \sum_{0 < |k| \le K} (\delta + K^{-1}) \left|\sum_{Q< q \le 2Q}\sum_{H_- \le h \le H_+} \int_{\eta q}^{\xi q} e \left( kqf(\alpha /q) - h\alpha \right) \, d\alpha \right|.$$

\noindent Since $H \ll  |k| \le K$, the error term in (\ref{N1 to N2}) satisfies
\begin{equation} \label{N1 to N2 error}\sum_{0 < |k| \le K} (\delta + K^{-1})\sum_{Q< q \le 2Q} \log(2+H) \ll  (\delta+K^{-1}) KQ \log K.
\end{equation}

By a change of variables, the integral in the expression for $N_2$ can be written as
$$ q \int_{\eta}^{\xi} e \left( q(kf(\beta) - h\beta) \right) \, d\beta.$$ 
The function $g(\beta) = q (k f(\beta) - h\beta)$ has second derivative $qkf''(\beta)$, which has modulus lying between constant multiples of $q|k|$.   Thus, by Lemma 4.4 of \cite{Titch86}, for any subinterval $\mathcal{I}$ of $[\eta, \xi]$, 
\begin{equation}  \int_\mathcal{I} e \left( q(kf(\beta) - h\beta) \right) \, d\beta \ll  \frac{1}{\sqrt{q|k|}}.
\label{integral est}\end{equation} 
Thus the contribution to $N_2$ from any $h$ with $H_- \le h \le h_-$ or $h_+ \le  h \le H_+$ is 
$$
\ll  \sum_{0<|k|\le K} (\delta + K^{-1})\sum_{Q< q \le 2Q} q \frac{1}{\sqrt{q|k|}} 
\ll  \delta K^{\frac{1}{2}} Q^{\frac{3}{2}} +K^{-\frac{1}{2}} Q^{\frac{3}{2}},
$$
and so we have 
\begin{equation}\label{N2 to N3 error}
N_2 = N_3  + O(\delta K^{\frac{1}{2}} Q^{\frac{3}{2}} +K^{-\frac{1}{2}} Q^{\frac{3}{2}}),
\end{equation}
where 
$$N_3  = \sum_{0 < |k| \le K}(\delta + K^{-1})\left|\sum_{Q< q \le 2Q} q \sum_{h_- < h < h_+} \int_{\eta }^{\xi } e \left( q(kf(\beta) - h\beta) \right) \, d\beta\right|.$$

Since $f'$ is continuous and $\inf k f'(\beta) < h_- < h < h_+ < \sup k f'(\beta)$, and since $f''$ is continuous and nonzero, it follows by the intermediate value theorem that there is a unique ${\beta_h = \beta_{k,h} \in [\eta, \xi]}$ such that $k f'(\beta_h) = h$.  Let
$$ \lambda_h = \lambda_{k,h}  = || k f(\beta_h) - h\beta_h ||. $$

\noindent By (\ref{integral est}), the terms of $N_3$ with $\lambda_h \le Q^{-1}$ contribute
$$ \ll  (\delta + K^{-1}) \sum_{0<|k|\le K}\  \sum_{\substack{ h_- < h < h_+ \\ \lambda_h \le Q^{-1}}} \  \sum_{Q < q \le 2Q} q^{\frac{1}{2}} |k|^{-\frac{1}{2}}.$$

\noindent By Lemma 2.3 of \cite{VaVe06} this is
$$  \ll  (\delta + K^{-1}) Q^{\frac{3}{2}} (K^{\frac{3}{2}} Q^{\varepsilon - 1} + K^{\frac{1}{2}} \log K),
$$
where $\varepsilon>0$ is any positive real number.  Thus we have 
\begin{equation}\label{N3 to N4 error}
N_3 = N_4 + O\left((\delta + K^{-1}) Q^{\frac{3}{2}} (K^{\frac{3}{2}} Q^{\varepsilon - 1} + K^{\frac{1}{2}} \log K\right),
\end{equation}
where $$N_4 =   \sum_{0 < |k| \le K}  (\delta + K^{-1})\left|\sum_{Q< q \le 2Q} q \sum_{\substack{h_- < h < h_+ \\ \lambda_h > Q^{-1}}} \int_{\eta }^{\xi } e \left( q(kf(\beta) - h\beta) \right) \, d\beta\right|.$$

Let $\beta_h = \beta_{k,h}$ be as above and let $\mu = (\xi - \eta)/2$.  Define
$$\mathcal{A}_1 := [\eta, \xi] \setminus [\beta_h - \mu, \beta_h + \mu] ,$$ 
$$\mathcal{A}_2 :=[\beta_h - \mu, \beta_h + \mu] \setminus [\eta, \xi] .$$
From the proof of Theorem 3 in \cite{VaVe06}, we see that for $i = 1,2$,
$$\int_{\mathcal{A}_i } e \left( q(kf(\beta) - h\beta) \right) \, d\beta \ll  \frac{1}{q(h-h_-)} + \frac{1}{q(h_+ -h)} .$$
Therefore we have
\begin{equation}\label{4to5}
N_4 = N_5 + O\left( (\delta + K^{-1})Q \sum_{0 < |k| \le K} \sum_{h_- < h < h_+} \frac{1}{(h-h_-)} + \frac{1}{(h_+ -h)} \right),
\end{equation}
where
$$N_5 = \sum_{0 < |k| \le K}(\delta + K^{-1})\left| \sum_{\substack{h_- < h < h_+ \\ \lambda_h > Q^{-1}}} \sum_{Q< q \le 2Q} q \int_{\beta_h -\mu}^{\beta_h+\mu } e \left( q(kf(\beta) - h\beta) \right) \, d\beta\right|.$$

\noindent Note that the error term in (\ref{4to5}) is 
\begin{equation} \label{N4 to N5 error}
 \ll  (\delta + K^{-1})Q \sum_{0 < |k| \le K}\log K  \ll  (\delta + K^{-1})QK\log K.
\end{equation}

We are left to deal with $N_5$.   Again following from the proof of Theorem 3 in  \cite{VaVe06}, we have that
$$ \sum_{Q< q \le 2Q} q \int_{\beta_h -\mu}^{\beta_h+\mu } e \left( q(kf(\beta) - h\beta) \right) \, d\beta \ll  Q^{1/2}\lambda_h^{-1}|k|^{-1/2} + Q^{(3-\theta)/2} |k|^{(-1-\theta)/2}.$$
Using Lemma 2.3 of  \cite{VaVe06} it then follows that
\begin{align}
N_5 & \ll  (\delta + K^{-1}) \sum_{0 < |k| \le K} \sum_{\substack{h_- < h < h_+ \\ \lambda_h > Q^{-1}}} \left( Q^{1/2}\lambda_h^{-1}|k|^{-1/2} + Q^{(3-\theta)/2} |k|^{(-1-\theta)/2} \right)\nonumber \\ &
\ll   (\delta + K^{-1}) \left( Q^{1/2+\varepsilon}K^{3/2} + Q^{3/2}K^{1/2}\log K + Q^{(3-\theta)/2} K^{(3-\theta)/2} \right)  \label{N5 error}.
\end{align}
\newline

We now have our upper bound for $N_1$.  Combining the error terms in (\ref{N1 to N2 error}), (\ref{N2 to N3 error}), (\ref{N3 to N4 error}), (\ref{N4 to N5 error}), and (\ref{N5 error}), we see that 
\begin{equation} \label{N1 bound}
 N_1 \ll  (\delta K + 1)Q\left(\frac{Q^{1/2}\log K }{K^{1/2}} +  \log K + \frac{K^{1/2}}{Q^{1/2-\varepsilon}} + (KQ)^{\frac{1-\theta}{2}}\right).
\end{equation}
Thus we see that
\begin{equation}\label{E messy}E(Q,\delta) \ll \frac{Q^2}{K} + (\delta K + 1)Q\left(\frac{Q^{1/2}\log K }{K^{1/2}} +  \log K + \frac{K^{1/2}}{Q^{1/2-\varepsilon}} + (KQ)^{\frac{1-\theta}{2}}\right).\end{equation}

The goal now is to find the choice of $K$  that minimizes $E(Q,\delta)$.  To simplify the computations, we allow $K\in \mathbb{R}$ for the time being.  We will take the floor function of our choice later to get back to $K\in \mathbb{N}$.  If $K > Q^{1-\frac{2}{3}\varepsilon}$, then
$$\delta K Q \left(\frac{K}{Q}\right)^{1/2}Q^{\varepsilon} > \delta Q^2,$$ and hence is too big to give an asymptotic formula.   Thus we may suppose that $K\le Q^{1-\frac{2}{3}\varepsilon}$.  Then, since $\theta <1$, we obtain 
\begin{equation}\label{E cleaner} E(Q,\delta) \ll  K^{-1} Q^2 + (\delta K + 1)Q\left(\left(\frac{Q}{K}\right)^{1/2}\log K + (KQ)^{(1-\theta)/2}\right).
\end{equation} 

If $\delta K \le 1$ then $K^{-1}Q^2 \ge \delta Q^2$, and we do not get our asymptotic formula.  So we assume that $\delta K > 1$ and (\ref{E cleaner}) simplifies to
\begin{equation}\label{E polished} E(Q,\delta) \ll  K^{-1} Q^2 + \delta K ^{1/2}Q^{3/2}\log Q + \delta(KQ)^{(3-\theta)/2}.
\end{equation} 
We replaced $\log K$ by $\log Q$ in the above bound to simplify our computations.  This is valid because the restrictions we have placed on $\delta$ and $K$ so far require that $\log K \ll  \log Q$.
The optimal choice for $K$ will occur when two of the three terms in (\ref{E polished}) are equal.  So we may reduce our analysis to three cases:  
$K = \delta^{-2/3}Q^{1/3}(\log Q)^{-2/3}, \ \ K = \delta^{\frac{-2}{5-\theta}}  Q^{\frac{1+\theta}{5-\theta}} , \ \ $ and $K = Q^{\frac{\theta}{2-\theta}}(\log Q)^{\frac{\theta}{2-\theta}}$.  These cases will yield three upper bounds for $E(Q,\delta)$.   We will then compare those bounds to find the least upper bound.
\newline
\newline
\emph{Case 1:} $K = \delta^{-2/3}Q^{1/3}(\log Q)^{-2/3}$.\newline
With this choice of $K$ we have that
$$K^{-1} Q^2 =  \delta K ^{1/2}Q^{3/2}\log Q= \delta^{2/3}Q^{5/3}(\log Q)^{2/3}$$
 and
$$\delta(KQ)^{(3-\theta)/2}  = \delta^{\theta/3}Q^{\frac{2}{3}(3-\theta)}(\log Q)^{-\frac{1}{3}(3-\theta)}.$$
Straightforward computations to find the dominating terms show that 
$$E(Q,\delta) \ll  \left\{ \begin{array}{l c l}
\delta^{2/3}Q^{5/3}(\log Q)^{2/3} & \mbox{ if } & \delta>>Q^{\frac{1-2\theta}{2-\theta}}(\log Q)^{-\frac{5-\theta}{2-\theta}} \\ 
\delta^{\theta/3}Q^{\frac{2}{3}(3-\theta)}(\log Q)^{-\frac{1}{3}(3-\theta)}& \mbox{ if } & \delta\ll Q^{\frac{1-2\theta}{2-\theta}}(\log Q)^{-\frac{5-\theta}{2-\theta}}.
\end{array}\right.$$
\newline
\newline
\emph{Case 2:} $K = \delta^{\frac{-2}{5-\theta}}  Q^{\frac{1+\theta}{5-\theta}}$.\newline
In this case we have
$$K^{-1} Q^2 = \delta(KQ)^{(3-\theta)/2}  = \delta^{\frac{2}{5-\theta}}Q^{\frac{3(3-\theta)}{5-\theta}}$$ and
$$\delta K ^{1/2}Q^{3/2}\log Q  = \delta^{\frac{4-\theta}{5-\theta}}Q^{\frac{8-\theta}{5-\theta}}(\log Q).$$
Thus, we see that 
$$E(Q,\delta) \ll  \left\{ \begin{array}{l c l}
\delta^{\frac{4-\theta}{5-\theta}}Q^{\frac{8-\theta}{5-\theta}}(\log Q) & \mbox{ if } & \delta>>Q^{\frac{1-2\theta}{2-\theta}}(\log Q)^{-\frac{5-\theta}{2-\theta}} \\ 
\delta^{\frac{2}{5-\theta}}Q^{\frac{3(3-\theta)}{5-\theta}}& \mbox{ if } & \delta\ll Q^{\frac{1-2\theta}{2-\theta}}(\log Q)^{-\frac{5-\theta}{2-\theta}}.
\end{array}\right.$$
\newline
\newline
\emph{Case 3:} $K = Q^{\frac{\theta}{2-\theta}}(\log Q)^{\frac{\theta}{2-\theta}}$.\newline
We now have 
$$\delta K^{1/2} Q^{3/2}\log Q = \delta (KQ)^{(3-\theta)/2} = \delta Q^{\frac{3-2\theta}{2-\theta}}(\log Q)^{\frac{4-3\theta}{2(2-\theta)}}$$ and
$$ K^{-1}Q^2 = Q^{\frac{4-\theta}{2-\theta}}(\log Q)^{\frac{\theta}{2(2-\theta)}}.$$
We obtain
$$E(Q,\delta) \ll  Q^{\frac{4-\theta}{2-\theta}}(\log Q)^{\frac{\theta}{2(2-\theta)}}.$$
\newline

Comparing the bounds from each of the three cases, we find that the least upper bound is given by
$$E(Q,\delta) \ll  \left\{ \begin{array}{l c l}
\delta^{2/3}Q^{5/3}(\log Q)^{2/3} & \mbox{ if } & \delta>>Q^{\frac{1-2\theta}{2-\theta}}(\log Q)^{-\frac{5-\theta}{2-\theta}} \\ 
\delta^{\frac{2}{5-\theta}}Q^{\frac{3(3-\theta)}{5-\theta}}& \mbox{ if } & \delta\ll Q^{\frac{1-2\theta}{2-\theta}}(\log Q)^{-\frac{5-\theta}{2-\theta}}.
\end{array}\right.$$
Hence we will choose $K = \lfloor\delta^{-\frac23}Q^{\frac13}(\log Q)^{-\frac23}\rfloor$ when $\delta>>Q^{\frac{1-2\theta}{2-\theta}}(\log Q)^{-\frac{5-\theta}{2-\theta}}$ and $K =  \lfloor\delta^{\frac{-2}{5-\theta}}  Q^{\frac{1+\theta}{5-\theta}} \rfloor$ when  $\delta\ll Q^{\frac{1-2\theta}{2-\theta}}(\log Q)^{-\frac{5-\theta}{2-\theta}}$.  Since we have an additional assumption that $\delta<1/2$, the first range for delta will only occur if $\theta>1/2$. This completes the proof of the theorem.  \hfill $\square$

\section{Proof of Theorem 2}

We obtain $N(Q,\delta)$ from $\widetilde{N}(Q,\delta)$ by a dyadic sum.  That is,
$$N(Q,\delta) = \sum_{r=1}^\infty \widetilde{N}(\frac{Q}{2^r},\delta).$$ 
It is easy to see that this sum converges since $\widetilde{N}(\frac{Q}{2^r},\delta) = 0$ if $2^{r-1}>Q$.  To avoid restrictions on $\delta$ in terms of $Q/2^r$, we will use the estimate for $E(Q,\delta)$ given by (\ref{E messy}).  We have
\begin{align*}
N(Q,\delta) & = \sum_{r=1}^\infty \widetilde{N}(\frac{Q}{2^r},\delta)
=  \sum_{r=1}^\infty \left(3 (\xi - \eta)\delta \left(\frac{Q}{2^r}\right)^2 +E(\frac{Q}{2^r},\delta) \right)
\\ & = \sum_{r=1}^\infty 3 (\xi - \eta)\delta \frac{Q^2}{4^r} + \sum_{r=1}^\infty F_r(Q,\delta),
\end{align*}
where 
$$F_r(Q,\delta) \ll  \frac{Q^2}{4^rK} + (\delta K + 1)\left(\frac{Q^{\frac32}\log K}{2^{\frac{3r}{2}}K^{\frac12}} + \frac{Q\log K}{2^r} + \frac{K^{\frac12}Q^{\frac12 + \varepsilon}}{2^{r(\frac12 + \varepsilon)}} + \frac{K^{\frac{1-\theta}{2}}Q^{\frac{3-\theta}{2}}}{2^{\frac{r(3-\theta)}{2}}}
\right).$$
Since $r$ only appears as an exponent of $(1/2)^\alpha$ for various values of $\alpha >0$, it is clear by the convergence of the geometric series that
\begin{align*}F(Q,\delta)& := \sum_{r=1}^\infty F_r(Q,\delta) \\ &\ll  \frac{Q^2}{K} + (\delta K + 1)Q\left(\frac{Q^{1/2}\log K }{K^{1/2}} +  \log K + \frac{K^{1/2}}{Q^{1/2-\varepsilon}} + (KQ)^{\frac{1-\theta}{2}}\right).\end{align*}
Note that this is the same estimate that is given for $E(Q,\delta)$ in (\ref{E messy}).  Thus the proof of Theorem 1 gives the bound for $F(Q,\delta)$.  We now return our attention to the main term of $N(Q,\delta)$.  We have
\begin{align*}
N(Q,\delta) & =  \sum_{r=1}^\infty 3 (\xi - \eta)\delta \frac{Q^2}{4^r} +F(Q, \delta) 
= 3 (\xi - \eta)\delta Q^2\frac{1/4}{1-1/4} + F(Q, \delta) \\ &
= (\xi - \eta)\delta Q^2 + F(Q, \delta),
\end{align*}
as desired. \hfill $\square$

\section{Proof of the corollaries}
 
Denote the piecewise upper bound given in (\ref{official error}) by $E_1(Q,\delta)$. To prove both corollaries, it is clearly enough to show that 
$$\frac{E_1(Q,\delta)}{\delta Q^2} \rightarrow 0$$
 as $Q\rightarrow\infty$ and $\delta\rightarrow0$.  We will call upon the assumption that $\delta \ge Q^{-\frac{1+\theta}{3-\theta}+\varepsilon}$.  When  $\delta>>Q^{\frac{1-2\theta}{2-\theta}}(\log Q)^{-\frac{5-\theta}{2-\theta}}$, we have 
$$\frac{E_1(Q,\delta)}{\delta Q^2} \ll  \frac{\delta^{2/3}Q^{5/3}(\log Q)^{2/3}}{\delta Q^2} = (\delta Q)^{-1/3}(\log Q)^{2/3} \le Q^{-\frac{2-2\theta}{3(3-\theta)}-\frac{\varepsilon}{3}}(\log Q)^{\frac{2}{3}},$$
which tends to $0$ as $Q\rightarrow\infty$.  Meanwhile, when $\delta\ll Q^{\frac{1-2\theta}{2-\theta}}(\log Q)^{-\frac{5-\theta}{2-\theta}}$, we have
$$\frac{E_1(Q,\delta)}{\delta Q^2} \ll  \frac{\delta^{\frac{2}{5-\theta}}Q^{\frac{3(3-\theta)}{5-\theta}}}{\delta Q^2} = \delta^{-\frac{3-\theta}{5-\theta}}Q^{-\frac{1+\theta}{5-\theta}} \le Q^{-\varepsilon\frac{3-\theta}{5-\theta}},$$
which also tends to $0$ as $Q\rightarrow\infty$.
\hfill $\square$

\subsection*{Acknowledgements}
The author would like to thank her adviser, Robert C. Vaughan, for his guidance throughout this research.

\end{document}